\documentclass[invmat]{svjour}
\usepackage{latexsym}
\usepackage{graphics}

\usepackage{graphicx}
\usepackage{amssymb}
\usepackage{amsfonts}
\usepackage[all]{xy}
\usepackage{amsmath}

\newcommand{\NN}{\mathbb{N}}

\newcommand{\al}{\alpha}

\newcommand{\A}{\mathcal{A}}

\newcommand{\C}{\mathcal{C}}

\newcommand{\R}{\mathcal{R}}

\begin{document}

\title{Symmetric product as moduli space of linear representations}

\author{Francesco Vaccarino\inst{1}
\thanks{The author is supported by research grant Politecnico di Torino n.119 - 2004, MSC 2000: 14L30, 13A50, 20G05}%
}

\institute{Politecnico di Torino - DSPEA, I-10129, Torino, Italy,
francesco.vaccarino@polito.it}

\date{26 august 2006}

\maketitle

\begin{abstract}
We show that the $n-$th symmetric product of an affine scheme
$X=\mathrm{Spec}\,A$ over a characteristic zero field is isomorphic as a scheme to the quotient by the general linear group of the scheme parameterizing $n-$dimensional linear representations of $A$. As a consequence we give generators and relations of the related rings of invariants as well as the equations of any symmetric products in term of traces. In positive characteristic we prove an analogous result for the associated varieties.\end{abstract}
%
\section{Introduction}\label{intro}
Let $F$ be an infinite field and let $F[y_1,\dots,y_n]^{S_n}$ be the
ring of \emph{symmetric polynomials} in $n$ variables. The general
linear group $GL(n,F)$ acts by conjugation on the full ring
$Mat(n,F)$ of $n\times n$ matrices over $F$. Denote by
$F[Mat(n,F)]^{GL(n,F)}$ the ring of the polynomial invariants for
this actions. It is well known that
\begin{equation}F[Mat(n,F)]^{GL(n,F)}\cong
F[y_1,\dots,y_n]^{S_n}.\label{cla}\end{equation} The above result
can be restated by saying that the scheme parameterizing
$n-$dimensional linear representations of $F[x]$ up to basis change
is isomorphic to the symmetric product $(\mathbb{A}^1)^n/S_n$ which
is defined as $\mathrm{Spec}\,F[y_1,\dots,y_n]^{S_n}$.

The main result of this paper is to show that the same thing happens
for any commutative algebra $A$ over a characteristic zero field.

\noindent Namely, when $F$ is a characteristic zero field, we prove
that the scheme $\R_n(A)//Gl(n,F)$ that parameterizes as a coarse
moduli space $n-$dimensional linear representations of a commutative
$F-$algebra $A$ up to basis change is isomorphic to the symmetric
product of its prime spectrum, i.e.
\[\R_n(A)//Gl(n,F)\cong
X^n/S_n=\mathrm{Spec}\,(A^{\otimes n})^{S_n},\] where
$X=\mathrm{Spec}\,A$.

Suppose $A$ is generated by say $m$ elements, then on geometric
points one can identify the symmetric product with a subscheme of
the $m-$tuples of diagonal $n\times n$ matrices and the above result
can be read as theorem of simultaneous diagonalization.

In positive characteristic we prove an analogous result for the associated varieties.

The proof is based on the carefully analysis of the morphism induced
by the composition $det\cdot \rho$ of the determinant with a
representation. This will be developed for flat algebra over a
commutative ring.

As a complement we give a presentation by generators and relations
of $(A^{\otimes n})^{S_n}$ that holds for $A$ flat and over any
commutative base ring $F$ giving then the equations of $X^{(n)}$ in
the flat case.

Thank to this we describe $V_n(A)^G\cong (A^{\otimes n})^{S_n}$ in
the characteristic zero case in term of traces and polynomial
identities and in positive characteristic in term of the coefficients of
the characteristic polynomial of polynomial of generic matrices.

These results generalizes the one given in {\cite{vf,v3}}.

\centerline{\textbf{Acknowledgements}}
I would like to thank M.Brion, C. De Concini and last but not least C.Procesi for useful discussions.

\section{Notation} Unless otherwise stated we adopt the following
notations:
\begin{itemize}
\item $F$ is a fixed base ring
\item we write {\it{algebra}} to mean {\it {commutative $F-$algebra}}
\item we denote by $\mathcal{C}_F$ the category of commutative $F-$algebras
\item $Sets$ the category of sets
\item we write $\A(B,C):=Hom_{\A}(B,C)$ in a category $\A$ with
$B,C\in Ob(\A)$ objects in $\A$.
\item for $A$ a set and any
additive monoid $M$, we denote by $M^{(A)}$ the set of functions
$f:A\rightarrow M$ with finite support.
\item let $\al\in M^{(A)}$,
we denote by $\mid \al \mid$ the (finite) sum $\sum_{a\in A}
\al(a)$,
\item given a set $I$ we denote by $\sharp I$ its cardinality.
\end{itemize}

\section{Representations}\label{rep}
Given an algebra $B$ write $Mat(n,B)$ for the ring of $n\times n$
matrices with entries in $B$.
\begin{definition}\label{rep}
 For a $n-$dimensional representation
of $A$ over $B$ one means an algebra homomorphism $A\to Mat(n,B)$.
We denote by $\R_n(A,B)$ the set of these representations. Given any
algebras $A,B$ and a homomorphism $\rho:A\to B$ we write $(\rho)_n$
for the induced map $Mat(n,A)\to Mat(n,B)$.
\end{definition}
The assignment $B\to \R_n(A,B)$ gives a covariant functor
$\R_n^A:\C_R\to Sets$ as can be easily checked.
\begin{proposition}({\cite{dp}} Sec.1\,)
For all $n\in\NN$ and $A\in\C_F$ there exist a unique algebra
$V_n(A)$ and a unique representation $\pi_n^A\in\R_n(A,V_n(A))$ such
that the map $\rho\mapsto(\rho)_n\cdot \pi_n^A$ is an isomorphism
\[\C_F(V_n(A),B)\xrightarrow{\cong}\R_n(A,B),\]
for all algebra $B$, i.e. the functor from schemes to sets
associated to $\R_n^A$ is represented by the affine scheme
$\R_n(A):=\mathrm{Spec}\,V_n(A)$
\end{proposition}
\begin{proof} Let $\Omega$ be a set and consider
$F[\xi_{ij,\omega}]$, a polynomial ring where $i,j=1,\dots,n$ and
$\omega\in\Omega$. Note that $F[\xi_{ij,\omega}]$ is isomorphic to
the symmetric algebra on the dual of $Mat(n,F)^{\Omega}$.

Let $A_{\Omega}=F\{x_{\omega}\}_{\omega\in\Omega}$ be the free
associative algebra on $\Omega$ then
\[\R_n(A_{\Omega},S)\cong Mat(n,S)^{\Omega} \cong
\C_F(F[\xi_{ij,\omega}],S)\] for any $S\in\C_F$. More precisely write
$D:=Mat(n,F[\xi_{ij,\omega}])$ and let $\xi_{\omega}\in D$ be given
by $(\xi_{\omega})_{ij}=\xi_{ij,\omega}$\,, $\forall\, i,j,\omega$,
these are called the ($n\times n$) {\em generic matrices} and were
introduced by C.\,Procesi (see {\cite{prolin}}).

Let $\pi:A_{\Omega}\to D$ be the $n-$dimensional representation
given by $x_{\omega}\mapsto \xi_{\omega}$. We have then that given
any $\rho\in \R_n(A_{\Omega},S)$, with $S\in\C_F$ there is a unique
$\bar{\rho}\in \C_F(F[\xi_{ij,\omega}],S)$ given by
$\xi_{ij,h}\mapsto (\rho(x_{\omega}))_{ij}$ and it is such that the
following diagram commutes
\begin{equation}
\xymatrix{
  A_{\Omega} \ar[dr]_{\rho} \ar[r]^{\pi}
                & B \ar[d]^{(\bar{\rho})_n}  \\
                & Mat(n,S)             }
\end{equation}
We now substitute $A_{\Omega}$ with an algebra $A\in\C_F$ and we let
\[0\longrightarrow J \longrightarrow A_{\Omega} \longrightarrow A
\longrightarrow 0\] be a presentation by generators and relations.

Let $I$ be the unique ideal in $F[\xi_{ij,\omega}]$ such that
\begin{equation}\label{ideal}
Mat(n,I)=D\,\pi(J)\,D
\end{equation}
then as one can easily check that, for all $S\in\C_F$
\begin{equation}
\R_n(R,S)\cong \C_F(F[\xi_{ij,\omega}]/I,S)
\end{equation}
via the lifting to a representation of $F$. Set now
$V_n(A):=F[\xi_{ij,\omega}]/I$. Let
$\{a_{\omega}\}_{\omega\in\Omega}$ be a set of generators of $A$, in
the same way as above we have a universal representation
\begin{equation}
\pi_n^A:\begin{cases}
        A\longrightarrow Mat(n,V_n(A)) \\
        r_{\omega}\mapsto \xi_{\omega}^A
      \end{cases}
\end{equation}
where $\xi_{\omega}^A$ is the image of $\xi_{\omega}$ via the
surjection $Mat(n,F[\xi_{ij,\omega}])\to Mat(n,V_n(A))$.
\end{proof}
\begin{remark} Note that $\R_n(A)$ could be quite complicated, as an example,
when $A=\mathbb{C}[x,y]$ we obtain that $\R_n(A)$ is the
\textit{commuting scheme} and it is not even known (but
conjecturally true) if it is reduced or not, see \cite{v3}.
\end{remark}
\begin{definition}
We set $\mathcal{G}_n(A):=\pi_n^A(A)\subset Mat(n,V_n(A))$.
\end{definition}
\begin{proposition}
The universal representation $\pi_n^A$ gives an isomorphism
$A\cong\mathcal{G}_n(A)$.
\end{proposition}
\begin{proof}
The representation $A\to Mat(n,A)$ given by $a\mapsto
diag(a,\dots,a)$ is injective so is $\pi_n^A$ by universality.
\end{proof}
\begin{definition}\label{coef}
The subalgebra of $V_n(A)$ generated by the coefficients of the
characteristic polynomial of the elements of $\mathcal{G}_n(A)$ will
be denoted by $C_n(A)$.
\end{definition}
\section{Symmetric Products}\label{sym}
Let $A$ be an algebra and $X=\mathrm{Spec}\,A$ its prime spectrum.
The symmetric group $S_n$ acts on the $n-$th tensor power
$A^{\otimes n}$ and as usual we write $TS^n_F(A)$ or simply
$TS^n(A)$ to denote the subalgebra of the invariants for this
actions, i.e. the symmetric tensors of order $n$ over $A$. The
$n-$th symmetric product of the affine scheme $X$ is the quotient
scheme of $X^n$ with respect to the above action and is usually
denoted by $X^{(n)}$. By definition
$X^{(n)}:=X^n/S_n:=\mathrm{Spec}\, TS^n(A)$.
\subsection{Polynomial
Laws}\label{poly} To link symmetric tensors to linear
representations we shall use the determinant so that we are lead to
the topic of polynomial laws: we recall the definition of this kind
of map between $F-$modules that generalizes the usual polynomial
mapping between free $F-$modules (see {\cite{bo,rl,rs}}).
\begin{definition} Let $A$ and $B$ be two
$F$-modules. A \emph{polynomial law} $\varphi$ from $A$ to $B$ is a
family of mappings $\varphi_{_{L}}:L\otimes_{F} A \longrightarrow
L\otimes_{F} B$, with $L\in\C_{F}$ such that the following diagram
commutes
\begin{equation}
\xymatrix{
  L\otimes_{F}A \ar[d]_{f\otimes id_A} \ar[r]^{\varphi_L}
                & L\otimes_{F} B \ar[d]^{f\otimes id_B}  \\
  M\otimes_{F} A \ar[r]_{\varphi_M}
                & M\otimes_{F} B             }
\end{equation}
for all $L,\,M\in\C(F)$ and all $f\in \C_{F}(L,M)$.
\end{definition}
\begin{definition}
Let $n\in \NN$, if $\varphi_L(au)=a^n\varphi_L(u)$, for all $a\in
L$, $u\in L\otimes_{F} A$ and all $L\in\C_{F}$, then $\varphi$ will
be said \emph{homogeneous of degree} $n$.
\end{definition}
\begin{definition}
If $A$ and $B$ are two $F$-algebras  and
\[
\begin{cases} \varphi_L(xy)&=\varphi_L(x)\varphi_L(y)\\
         \varphi_L(1_{L\otimes A})&=1_{L\otimes B}
         \end{cases}
         \]
for $L\in\C_{F}$ and for all $x,y\in L\otimes_{F} A$, then $\varphi$
is called \emph{multiplicative}.
\end{definition}
\begin{remark} A polynomial law $\varphi:A\rightarrow B$  is a
natural transformation $-\otimes_{F}A\rightarrow -\otimes_{F}B$.
\end{remark}
Let $A$ and $B$ be two $F$-modules and $\varphi:A\rightarrow B$ be a
polynomial law. The following result on polynomial laws is a
restatement of Th\'eor\`eme I.1 of {\cite{rl}}.
\begin{theorem}\label{roby} Let $S$ be a set.
\begin{enumerate}
\item Let $L=F[x_s]_{s\in S}$ and let $\{a_{s}\, :\,s\in S\}\subset A$
be such that $a_{s}=0$ except for a finite number of $s\in S$, then
there exist $\varphi_{\xi}((a_{s}))\in B$, with $\xi \in \NN^{(S)}$,
such that:
\[\varphi_{_{L}}(\sum_{s\in S} x_s\otimes
a_{s})=\sum_{\xi \in \NN^{(S)}}  x^{\xi}\otimes
\varphi_{\xi}((a_{s})),\] where $x^{\xi}:=\prod_{s\in S}
x_s^{\xi_s}$.
\item Let $R$ be any commutative $F$-algebra and let
$(r_s)_{s\in S}\subset R$, then: \[\varphi_{_{R}}(\sum_{s\in S}
r_s\otimes a_{s})=\sum_{\xi \in \NN^{(S)}}  r^{\xi}\otimes
\varphi_{\xi}((a_{s})),\] where $r^{\xi}:=\prod_{s\in S}
r_s^{\xi_s}$.
\item If $\varphi$ is homogeneous of degree $n$, then one has
$\varphi_{\xi}((a_{s}))=0$ if $\mid \xi \mid$ is
different from $n$. That is: \[\varphi_{_{R}}(\sum_{a\in A}
r_a\otimes a)=\sum_{\xi \in \NN^{(A)},\,\mid \xi \mid=n}
r^{\xi}\otimes \varphi_{\xi}((a)).\] In particular, if $\varphi$ is
homogeneous of degree $0$ or $1$, then it is constant or linear,
respectively.
\end{enumerate}
\end{theorem}
\begin{remark}\label{coef}
The above theorem means that a polynomial law $\varphi:A\rightarrow
B$ is completely determined by its coefficients
$\varphi_{\xi}((a_{s}))$, with $(a_s)_{s\in S} \in S^{(\NN)}$.
\end{remark}
\begin{remark} If $A$ is a free $F$-module and $\{a_{t}\,
:\, t\in T\}$ is a basis of $A$, then $\varphi$ is completely
determined by its coefficients $\varphi_{\xi}((a_{t}))$, with $\xi
\in \NN^{(T)}$. If also $B$ is a free $F$-module with basis
$\{b_{u}\, :\, u\in U\}$, then $\varphi_{\xi}((a_{t}))=\sum_{u\in
U}\lambda_{u}(\xi)b_{u}$. Let $a=\sum_{t\in T}\mu_{t}a_{t}\in A$.
Since only a finite number of $\mu_{t}$ and $\lambda_{u}(\xi)$ are
different from zero, the following makes sense:
\begin{eqnarray*}\varphi(a)=\varphi(\sum_{t\in T}\mu_{t}a_{t}) =  \sum_{\xi\in
\NN^{(T)}} \mu^{\xi}\varphi_{\xi}((a_{t}))& = & \sum_{\xi\in
\NN^{(T)}} \mu^{\xi}(\sum_{u\in U}\lambda_{u}(\xi)b_{u})\\ & = &
\sum_{u\in U}(\sum_{\xi\in \NN^{(T)}}\lambda_{u}(\xi)
\mu^{\xi})b_{u}.\end{eqnarray*} Hence, if both $A$ and $B$ are free
$F$-modules, a polynomial law $\varphi:A\rightarrow B$ is simply a
polynomial map.
\end{remark}
\begin{definition}\label{polset}
Let $A,B\in\C_F$ be two algebras we write $P_n(A,B)$ for the
multiplicative homogeneous polynomial mapping $A\to B$ of degree
$n$.
\end{definition}
The assignment $B\to P_n(A,B)$ determines a functor $P_n^A:\C_F\to
Sets$ as one can easily check.

\subsection{Symmetric products as representing schemes}\label{sym1}
From now on $A$ will be a flat algebra. There is an element
$\gamma_n\in P_n(A,TS^n(A))$ given by $\gamma_n(a)= a^{\otimes n}$
such that the composition $\C_F(TS^n(A),B)\ni\rho\mapsto\rho\cdot\gamma_n\in P_n(A,B)$ gives an isomorphism
\begin{equation}\C_F(TS^n(A),B)\xrightarrow{\cong}P_n(A,B),\end{equation}
 i.e. the
functor $P_n^A$ is represented by the symmetric product $X^{(n)}$,
(see {\cite{bo}} chap.IV). We call $\gamma_n$ {\textit{the universal
mapping of degree n}}.

The proof goes as follows: by Rem.~\ref{coef} a polynomial law is
determined by its coefficients. If $A$ is free then one can extract
from the set of coefficients of $\gamma_n$ a linear basis of
$TS^n(A)$. Then any polynomial law homogeneous of degree $n$
correspond to a specialization of these coefficients of $\gamma_n$.
The requirement to be multiplicative corresponds to the requirement
the above specialization to be an algebra homomorphism, since
$\gamma_n$ is obviously multiplicative. Then one applies the good
properties of $TS^n$ with respect to inverse limits and get the
desired result for flat algebras, see \cite{laz}.

\subsection{Generators}
Let $a\in A$, there is an algebra homomorphism
$\eta_a:F[x_1,\dots,x_n]^{S_n}\cong TS^n(F[x])\to TS^n(A)$ induced
by the evaluation of $x$ at $a$. We write $e_i^n(a)=\eta_a(e_i^n)$
where $e_i^n$ is the $i-$th elementary symmetric polynomial in $n$
variables. Given and independent variable $t$ we have an induced map
$\overline{\eta}_a:F[t][x_1,\dots,x_n]^{S_n}\cong TS^n(F[t])\to
TS^n(A[t])$ such that
\begin{equation}
t^n+\sum_ie_i^n(a)=\overline{\eta}_a(t^n+\sum_ie_i^n)=\overline{\eta}_a(\prod_i(t+x_i))=(t+a)^{\otimes
n}
\end{equation}
so that $e_i^n(a)$ is the orbit sum of $a^{\otimes i}\otimes
1^{\otimes n-i}$. Note that $e_n^n(a)=a^{\otimes n}=\gamma_n(a)$.

\begin{proposition}[Generators]
Let $A$ be a commutative flat algebra generated by $\{a_i\}_{i\in
I}$ then $TS^n_F(A)$ is generated by $e_i^n(a)$ where $a=\prod_i
r_i^{\al_i}$ is such that $\sum_i\al_i\leq max(n,\sharp I(n-1))$ and
the $\al_i$ are coprime.
\end{proposition}\label{gen}
\begin{proof} In \cite{vf} Th.1 we prove the above statement for
$A=P$ a free polynomial $F-$algebra. Given $A=P/I$ flat then
$TS^n_F(P)\to TS^n_{F}(A)$ is onto and we are done.
\end{proof}

\begin{proposition}\label{newt}
Let $F\supset \mathbb{Q}$ with $\mathbb{Q}$ the rational integers
and let $A$ be generated by $\{a_i\}_{i\in I}$ then $TS^n_F(A)$ is
generated by $e_1^n(a)$ where $a=\prod_i r_i^{\al_i}$ is such that
$\sum_i\al_i\leq n$.
\end{proposition}
\begin{proof}
It follows from Prop.\ref{gen}, Newton formulas and Noether's bound.
\end{proof}

\section{Determinant and isomorphism}
We prove the key result of this paper.
\begin{theorem}\label{main}
The composition $det\cdot \pi_n^A$ induces an isomorphism\\
$$\delta_n^A:TS^n_F(A)\to C_n(A)$$ for all flat $A\in \C_F$.
\end{theorem}
\begin{proof}
By Sect.~\ref{sym1} there is a unique algebra homomorphism
$\delta_n^A:TS^n_F(A)\to C_n(A)\subset A_n(A)$ such that
$det\cdot\pi_n^A=\delta_n^A\cdot\gamma_n$. We have an identification
\[1_{F[t]}\otimes \delta_n^A(t\otimes 1_A+1\otimes a)^{\otimes n}=det(t\otimes
1_A+1\otimes\pi_n^A(a))\]
for all $a\in A$ with $t$ and independent
variable. By definition of $C_n(A)$ and Prop.\ref{gen} one has that $\delta_A^n$ is surjective
By Sect.~\ref{rep} the representation $\partial_n:A\to
Mat(n,A^{\otimes n})$ given by $\partial(a)= diag(a\otimes
1^{\otimes n-1},\,1\otimes a\otimes 1^{\otimes
n-2},\,\dots\,,\,1^{\otimes n-1}\otimes a)$ corresponds uniquely to
a homomorphism $\overline{\partial_n}:V_n(A)\to A^{\otimes n}$ such
that $\partial_n=(\overline{\partial_n})_n\cdot \pi_n^A$. Observe
that
$\overline{\partial_n}\cdot\det\cdot\pi_n^A=det\cdot\partial_n(a)=a^{\otimes
n}=\gamma_n(a)$ hence the restriction of $\partial_n$ to $C_n(A)$ is surjective by Prop.\ref{gen} and gives an inverse to $\delta_n^A$ and we are done.
\end{proof}

\section{Invariants}\label{inv}
We want to study in a GIT fashion the equivalence classes of the
representations of $A$ under basis changes, i.e. under the action of
the general linear group $G:=Gl(n,F)$. The right object is the
categorical quotient
$\mathcal{R}_n(A)//G:=\mathrm{Spec}\,(V_n(A)^G)$, where as usual
$V_n(A)^G$ denotes the invariants for the $G-$action induced on
$V_n(A)$ by basis change on $F^n$ .
\begin{theorem}
If $F$ is an algebraically closed field the $F-$points of the variety associated to $X^{(n)}$ are in one to one correspondence with the equivalence classes of the semisimple $n-$dimensional linear representations of $A$.
\end{theorem}
\begin{proof}
It follows from Th.\ref{main} plus Sec.4\cite{pis}.
\end{proof}
\begin{theorem}\label{ch0}
Let $F\supset \mathbb{Q}$ then
$$\delta_n^A:TS^n_F(A)\cong V_n(A)^G$$ i.e.
$$X^{(n)}\cong \mathcal{R}_n(A)//G$$
\end{theorem}
\begin{proof}
In \cite{v3} we proved that the statement is true when $A$ is a free
polynomial ring. Then the characteristic zero case follows by
Th.\ref{main} and the reductivity of $G$.

\noindent {\it{Another proof}.} The statement also follows observing
that, by C.\,Procesi \cite{p1} one has
$F[\xi_{ij,\omega}]^G=C_n(A_{\Omega})$ and again the result follows
by Th.\ref{main} and the reductivity of $G$.
\end{proof}
\begin{remark}
Let $A=F[x]$ then Th.\ref{ch0} implies $$F[x_1,\dots,x_m]^{S_n}\cong
F[Mat(n,F)]^G$$ a well know result.
\end{remark}
\begin{remark}
If $A$ is reduced then Th.\ref{ch0} implies that $V_n(A)^G$ is
reduced too. This gives some support to the conjecture that
$V_n(\mathbb{C}[x,y])$ is reduced.
\end{remark}

\subsection{Equations}
Let now $F$ be again an arbitrary commutative ring. Suppose you have
an homomorphism $f:A\to B$ of algebras then on can easily check that
the kernel of $TS^nf:TS^n(A)\to TS^n(B)$ is linearly generated by
the orbit sums (under $S_n$) of elements $a_1\otimes\cdots\otimes
a_n$ such that $\exists k\in \{1,\dots,n\}$ with $a_k\in \ker f$.
Now any such element can be expressed as a polynomial in the
$e_i^n(a)$ with $a$ varying into the set of monomials in the $a_j$
(see \cite{vf} Lemma 1.2 and Cor.2.3). Thus we have the following
\begin{proposition}\label{ker1}
Let $f:A\to B$ be an algebra homomorphism then the kernel of
$TS^nf:TS^n(A)\to TS^n(B)$ is generated as an ideal by the elements
$e_i^n(a)$ with $i=1,\dots,n$ and $a\in \ker f$.
\end{proposition}
\begin{corollary}
Suppose $F\supset \mathbb{Q}$ then $\ker TS^nf$ is generated by
$e_1^n(a)$ with $a\in \ker f$.
\end{corollary}
\begin{proof}
It follows from the above Proposition by Newton's formulas.
\end{proof}

Let $P_{\Omega}=F[x_{\omega}]_{\omega\in\Omega}$ be the free
polynomial algebra on $\Omega$. We set $P^+_{\Omega}$ for the
augmentation ideal i.e. the kernel of the evaluation
$x_{\omega}\mapsto 0$ for all $\omega$.

Consider the polynomial ring $F[e_{i,\mu}]$ freely generated by the
symbols $e_{i,\mu}$ with $1\geq i$ and $\mu$ that varies in the set
of monomials $\prod_{\omega}x_{\omega}^{\alpha_{\omega}}$ having
coprime exponents. By Prop.\ref{gen} for all $n$ there is a
surjective homomorphism $\kappa_n:F[e_{i,\mu}]\to TS^n(P_{\Omega})$
given by $e_{i,\mu}\mapsto e_i^n(\mu)$ if $i\leq n$ and
$e_{i,\mu}\mapsto 0$ for $i>n$.

Given $f\in P_{\Omega}^+$ we can compute $e_m^n(f)=$ orbit sum under
$S_n$ of $f^{\otimes m}\otimes 1^{\otimes n-m}$ for all $n\geq m$
(for $n<m$ it is zero) and express it as a polynomial in the
$e_i^n(\mu)$ with $\mu$ monomial with coprime exponents. It has been
shown that this expression stabilizes as $n>>m$ (see \cite{vf}
Prop.3.4 or \cite{vz} 9.1). Thus we have a well defined polynomial
law $\overline{e}_m:P_{\omega}^+\to F[e_{i,\mu}]$ homogeneous of
degree $m$. In \cite{vf} we prove that $\ker\kappa_n$ is linearly
generated by the coefficients of $\overline{e}_m$ with $m>n$.

Let us give an improvement of that result.
\begin{theorem}[Relations]\label{rel}
\begin{enumerate}
\item Let $A=P_{\omega}/I$ be a flat algebra, then the kernel of the surjection
$F[e_{i,\mu}]\xrightarrow{\kappa_n}TS^n(P_{\Omega})\to TS^n(A)$ is
linearly generated by the coefficients of $\overline{e}_m$ with
$m>n$ plus the lifting to $F[e_{i,\mu}]$ of $e_k^n(g)$ where $k\leq
n$ and $g\in I\subset P_{\Omega}$.
\item Suppose $F$ is an infinite field: the kernel of $\kappa_n$ is
generated as an ideal by $e_m(f)$ with $m>n$ and $f\in
P_{\Omega}^+$.
\item Suppose $F$ is an infinite field and let $A=P_{\omega}/I$.
The kernel of the surjection
$F[e_{i,\mu}]\xrightarrow{\kappa_n}TS^n(P_{\Omega})\to TS^n(A)$ is
generated as an ideal by $e_m(f)$ with $m>n$ and $f\in P_{\Omega}^+$
plus the lifting to $F[e_{i,\mu}]$ of the $e_k^n(g)$ where $k\leq n$
and $g\in I\subset P_{\Omega}$.
\end{enumerate}
\end{theorem}
\begin{proof}\,

\begin{enumerate}
\item It follows from the above discussion and Prop.\ref{ker1}
\item A linear form that is zero on the linear subspace
generated by $e_k(f)$ with and $f\in P^+_{\Omega}$ is zero also on
the subspace generated by the coefficients as $f$ varies.
\item By \textbf{1.} and \textbf{2.}
\end{enumerate}
\end{proof}
\begin{remark}
The above Theorem give the equations of the symmetric product of any
affine scheme.
\end{remark}

\subsection{Traces} Along this paragraph we suppose
again $F\supset \mathbb{Q}$. When we have a polynomial law $\varphi$
homogeneous of degree $n$ we can consider its full polarization
$\varphi_{\mathbf{1}_n}$ that is the coefficient of $t_1t_2\cdots
t_n$ in $\varphi(t_1 x_1+\dots +t_1 x_n)$ where $t_i$ are commuting
independent variables and the $x_j$ are generic elements in the
domain of $\varphi$. It is well known that the process of
polarization is effective via restitution
$n!\varphi(x)=\varphi_{\mathbf{1}_n}(x,x,\dots,x)$.

\noindent What happens for $e_n^n$? We write \[e_n^n(t_1 x_1+\dots
+t_n x_n)=t_1t_2\cdots t_n e_{\mathbf{1}_n}^n(x_1,\dots,x_n) + Z\]
and if one observe that $e_n^n(t_1 x_1+\dots +t_n x_n)=(t_1
x_1+\dots +t_n x_n)^{\otimes n}$ then it follows easily that
$e_{\mathbf{1}_n}^n(x_1,\dots,x_n)=\sum_{\sigma\in
S_n}x_{\sigma(1)}\otimes\cdots\otimes x_{\sigma(n)}$\,. We would
like now to express $e_{\mathbf{1}_n}^n(x_1,\dots,x_n)$ in terms of
$e_1^n(\mu)$ with $\mu$ a monomial in the $x_i$.

It is clear that $\delta_n^A(e_i^n(a))=\psi_i^n(\pi_n^A(a))$ where
\begin{equation}det(t+\pi_n^A(a))=t^n+\sum_i\psi_i^n(\pi_n^A)t^{n-i}\end{equation}
in particular $e_1^n(a)$ is identified with the trace of
$\pi_n^A(a)$ and $e_n^n(a)$ with the determinant. It is well known
that (see \cite{p1}) the full polarization $\chi_{\mathbf{1}_n}$ of
the determinant can be expressed as a special polynomial in
traces of monomials, namely consider the cycle decomposition
$\sigma=\sigma_1\cdots\sigma_k\in S_n$ and let correspond to it the
product $T_{\sigma}=tr(\mu_1)\cdots tr(\mu_k)$ where
$\mu_h=x_{h_1}\cdots x_{h_l}$ being $\sigma_h=(h_1\,h_2\,\dots
h_l)$, then
\begin{equation}
\Psi_{\mathbf{1}_n}=\sum_{\sigma\in S_n}\epsilon_{\sigma}T_{\sigma}
\end{equation}
and it is a well celebrated theorem due to Procesi \cite{p1} and
Razmyslov\cite{ra} that all the relations (in characteristic zero)
between the invariants of matrices, i.e. between traces of monomial
of generic $n\times n$ matrices are consequences, in the sense of
$T-$ideals, of $\Psi_{\mathbf{1}_{n+1}}$.

Let us summarize all we are able say in the characteristic zero
case.
\begin{theorem}
Let $A=P_{\Omega}/I$ be a commutative $F-$algebra. The ring of the
invariants $V_n(A)^G$ is generated by traces of monomial of generic
matrices $\xi^A$ and the ideal of relations is generated by the
evaluation of $\Psi_{\mathbf{1}_{n+1}}$ at the elements of
$\mathcal{G}_n(A)$ and by the traces $tr(\pi_n^A(f))$ with $f\in I$.
The same obviously holds \textit{mutatis mutandis} in $TS^n_F(A)$
\end{theorem}
\begin{proof}
It follows for the above discussion, Th.\ref{rel}.3 using Newton's
formulas.
\end{proof}

\section{Positive characteristic}
In this subsection $F$ will be an infinite field of arbitrary
characteristic. Set $N_n(R)$ for the nilradical of $V_n(A)$ and
$\R_{n,\,red}(A):= \mathrm{Spec}\,V_n(R)/N_n(R)$ for the reduced
scheme associates to $\R_n(A)$. Since the nilradical it is
$G$-stable then the action of $G$ on $\R_n(A)$ can be restricted to
the variety $\R_{n,\,red}(A)$.
\begin{theorem}
The above isomorphism gives
$$TS^n_R(A)_{red}\cong V_n(A)_{red}^G$$
i.e. $$X^{(n)}_{red}\cong \mathcal{R}_{n,red}(A)//G.$$
\end{theorem}
\begin{proof}
The injectivity follows passing to an algebraically closed field and
then observing that a tuple of commuting matrices can be put
simultaneously in upper triangular form. One can then reach a tuple
of diagonal matrices via a one parameter subgroup of $G$. Hence if
an invariant regular function is zero on a tuple of diagonal
matrices then $f=0$.
Surjectivity comes from $TS^n_F F[x_1,\dots,x_m]\twoheadrightarrow TS^n_FA$ and Cor.4.1 in {\cite{v3}}
that states that the morphism $F[Mat(n,F)]^G\to  TS^n_F F[x_1,\dots,x_m]$ induced by restriction to diagonal matrices is onto for any commutative ring $F$.
\end{proof}

\end{document}